# SOME RELATIONS BETWEEN MUTUAL INFORMATION AND ESTIMATION ERROR IN WIENER SPACE

BY EDDY MAYER-WOLF[1] AND MOSHE ZAKAI

*Technion*

The model considered is that of "signal plus white noise." Known connections between the noncausal filtering error and mutual information are combined with new ones involving the causal estimation error, in a general abstract setup. The results are shown to be invariant under a wide class of causality patterns; they are applied to the derivation of the causal estimation error of a Gaussian nonstationary filtering problem and to a multidimensional extension of the Yovits–Jackson formula.

**1. Introduction.** The classical "additive Gaussian channel" model consists of an $m$-dimensional "white noise" $\{n_t, t \in [0,T]\}$, an $m$-dimensional (not necessarily Gaussian) independent "signal process" $\{x_t, t \in [0,T]\}$ and the "received signal" $y_t = \sqrt{\gamma} x_t + n_t$, where $\gamma$ is the signal to noise parameter. (It also deals with the stationary version where $[0,T]$ is replaced by $(-\infty, \infty)$ and $x_t$ is assumed to be a stationary process.) In the context of filtering theory, the main entities are the noncausal estimate and its associated estimation mean square error

$$(1.1) \qquad \widetilde{\varepsilon}^2(\gamma) = \int_0^T E|x_t - E(x_t|y_\theta, \theta \in [0,T])|^2 \, dt$$

as well as the causal estimate and its associated filtering mean square error

$$(1.2) \qquad \widehat{\varepsilon}^2(\gamma) = \int_0^T E|x_t - E(x_t|y_\theta, \theta \in [0,t])|^2 \, dt.$$

Another aspect of the white Gaussian channel is the "mutual information" $I(x,y)$ between the signal process and the received message defined by

$$(1.3) \qquad I(x,y) = E \log \frac{dP(x,y)}{d(P(x) \times P(y))}$$

Received October 2006; revised January 2007.
[1]Supported in part by the Technion VPR Fund.
*AMS 2000 subject classifications.* Primary 60G35, 94A17; secondary 28C20, 60H35.
*Key words and phrases.* Estimation errors, mutual information, abstract Wiener space, resolution of identity, causality.

where the argument of the logarithm is the Radon–Nikodym derivative between the joint measure of $x.$ and $y.$ and the product measure induced by $x.$ and $y..$ This notion was introduced by Shannon and is essential in the definition of channel capacity, which in turn determines the possibility of transmitting signals through the channel with arbitrarily small error. The mutual information between "random objects" has been thoroughly analyzed and explicit results have been obtained, particularly for Gaussian signals and noise (cf. [7]).

Recently, Guo, Shamai and Verdu [3] derived interesting new results for the Gaussian channel relating the mutual information with the noncausal estimation error. These results were extended in [13] to include the abstract Wiener space setup, thus extending considerably the applicability of the new relations. As for the causal estimation problem, some general results are known, starting with the Yovits–Jackson formula [12], see Snyders [8, 9] for further results in this direction. Moreover, the relation between mutual information and the causal error appeared in the literature in the early 1970s [1, 5]. The possibility of extending these results to the abstract Wiener space was pointed out in [13].

The purpose of this paper is to consider the "noise" as a general Gaussian random vector and to establish connections between the causal estimation error and mutual information in this abstract setting. In addition, some new consequences of these connections are obtained, such as the concavity of the causal estimation error as a function of the noise-to-signal ratio (Corollary 3.3) as well as an explicit expression for the causal error in the estimation of a general (not necessarily stationary) Gaussian signal (Theorem 4.1), from which the Yovits–Jackson formula for a stationary Gaussian signal process follows quite directly (Proposition 4.3).

The context of an abstract Wiener space, apart from its intrinsic elegance, accommodates a wide range of signal models involving, for example, vector valued processes time reversed in some of its coordinates. We feel that this flexibility justifies the inclusion of the necessary abstract and sometimes tedious Wiener space analysis background in Section 2 and Section 3.1. On the other hand, as pointed out in the next section, the main results can also be of value to the reader who prefers to interpret their ingredients as concrete one dimensional processes.

We now outline the contents of this paper. In the next section the basic abstract Gaussian channel setup is introduced and some preliminary adaptedness results in the associated abstract Wiener space are established. In Section 3 the results of [1] and [5] are extended to the abstract Wiener space which, however, does not have any intrinsic notion of causality. Accordingly, it is equipped with a time structure by adding an appropriate "chain of projections" (namely, a continuous increasing resolution of the identity). It turns out that the causal estimation error is independent of the particular choice of the chain of projections, and is closely related to the mutual information $I(x,y)$. Moreover, this relation persists when the independence assumption between the signal $x$ and the noise $n$ is relaxed to allow for nonanticipative dependence, as in [5]. These results when combined with the earlier results on the nonadapted error yield a direct relation between the causal and noncausal errors. In Section 4 we derive the formulae alluded to in the previous paragraph, namely $\hat{\varepsilon}^2(\gamma) = \gamma^{-1} \sum_i \log(1 + \lambda_i \gamma)$ for a Gaussian process $x_t$ on $[0,T]$ whose correlation function has an eigenfunction expansion $\sum_i \lambda_i \varphi_i(s)\varphi_j(t)$, and the multidimensional version of the Yovits–Jackson formula $\hat{\varepsilon}^2(\gamma) = (2\pi\gamma)^{-1} \int_{-\infty}^{\infty} \log \det(I + \gamma\sigma(\xi))\, d\xi$ for a stationary Gaussian signal with (matricial) spectral density $\sigma$.

**2. Preliminaries.** This work studies the basic signal plus noise model, which will now be formally described, modeled on the abstract Wiener space to allow for maximal generality as mentioned in the Introduction. However, many of the paper's statements—including Theorem 3.1, Corollary 3.3 and the contents of Section 4—can be appreciated even in the simplest instance [cf. with (2.5)]

(2.1) $$y_t = u_t + w_t, \qquad 0 \leq t \leq T$$

(where the noise is represented by the Brownian motion $\{w_t\}$ and, at each $t \in [0,T]$, the signal $u_t$ depends at most on a "hidden" process $\{x_t\}$ independent of $\{w_t\}$ and, via feedback, on $y$'s "past" $\{y_s, 0 \leq s \leq t\}$, i.e., $u_t = U(x_0^T, y_0^t))$, without the need to master the details of the abstract setup whose data we now list:

M1. A complete filtered probability space $(\Theta, \mathcal{F}, \{\mathcal{F}_t, 0 \leq t \leq 1\}, P)$.

M2. A random variable $\mathbf{x}$ defined on $(\Theta, \mathcal{F}, P)$ taking its values in a Polish space $X$ and inducing on it its image measure $\mu_{\mathbf{x}}$.

M3. A centered nondegenerate Gaussian random variable $\mathbf{w}$ defined on $(\Theta, \mathcal{F}, P)$, independent of $\mathbf{x}$, taking values in a Banach space $\Omega$ with image measure $\mu_{\mathbf{w}}$, and separable associated reproducing kernel Hilbert space $H$. The non-degeneracy assumption means that $H$ is densely embedded in $\Omega$, namely, $(\Omega, H, \mu_{\mathbf{w}})$ is an abstract Wiener space and $\Omega^* \subset_{>} H \subset_{>} \Omega$.

M4. A time structure on $(\Omega, H, \mu_{\mathbf{w}})$ in the form of a continuous strictly increasing coherent resolution of the identity $\{\pi_t, 0 \leq t \leq 1\}$ of $H$, namely a (continuous, increasing) family of orthogonal projections on $H$ ranging from $\pi_0 = \mathbf{0}_H$ to $\pi_1 = \mathrm{Id}_H$, such that $\pi_t \Omega^* \subset \Omega^*$ and $_\Omega \langle \mathbf{w}, \pi_t l \rangle_{\Omega^*}$ is $\mathcal{F}_t$-adapted, for all $0 \leq t \leq 1$ and $l \in \Omega^*$.

With such a time structure one can mimic the standard resolution of identity $(\pi_t h). = h._{\wedge t}$ in classical Wiener space $C_0[0,1]$ (in fact (cf. [11], Theorem 5.1) any abstract Wiener space thus equipped with a resolution of the identity is equivalent "in a suitable sense" to $C_0([0,1];\mathbb{R}^d)$, for some $d \in \mathbb{N} \cup \infty$. This will not be used in the sequel):

(i) Any $\Omega$-valued random variable $\mathbf{z}$ induces a filtration $\{\mathcal{F}^{\mathbf{z}}_t, 0 \le t \le 1\}$ in $(\Theta, \mathcal{F})$

(2.2) $$\mathcal{F}^{\mathbf{z}}_t = \sigma(_\Omega\langle \mathbf{z}, \pi_t l\rangle_{\Omega^*}, l \in \Omega^*), \qquad 0 \le t \le 1.$$

(The above adaptedness requirement can be expressed as $\mathcal{F}^{\mathbf{w}}_\cdot \subset \mathcal{F}_\cdot$.)

(ii) Given a generic subfiltration $\{\mathcal{G}_t, 0 \le t \le 1\}$ of $\{\mathcal{F}_t, 0 \le t \le 1\}$, an $\Omega$-valued random variable $\mathbf{z}$ is said to be $(\pi_\cdot, \mathcal{G}_\cdot)$-*adapted* if $_\Omega\langle \mathbf{z}, \pi_t l\rangle_{\Omega^*}$ is $\mathcal{G}_t$-measurable, for all $0 \le t \le 1$ and $l \in \Omega^*$. Examples of $(\pi_\cdot, \mathcal{G}_\cdot)$-adapted random variables are provided in increasing generality, for a partition $\{0 = t_0 < \cdots < t_n = 1\}$ of $[0,1]$ (and with $\pi^t_s := \pi_t - \pi_s$) by

(2.3) $$\mathbf{h} = \sum_{k=0}^{n-1} a_k h_k, \qquad a_k \in L^2(\Theta, \mathcal{G}_{t_k}, P), h_k \in \pi^{t_{k+1}}_{t_k}(H),$$

(2.4) $$\mathbf{h} = \sum_{k=0}^{n-1} \mathbf{h}_k, \qquad \mathbf{h}_k \in L^2(\Theta, \mathcal{G}_{t_k}, P; \pi^{t_{k+1}}_{t_k}(H)),$$

(iii) A mapping $g: \Omega \to \Omega$ is $\pi_\cdot$-*nonanticipative* if $g(\mathbf{z})$ is $(\pi_\cdot, \mathcal{F}^{\mathbf{z}}_\cdot)$-adapted for any $\Omega$-valued $\mathbf{z}$, that is, if $_\Omega\langle g(\mathbf{z}), \pi_t l\rangle_{\Omega^*}$ is $\mathcal{F}^{\mathbf{z}}_t$-measurable for all such $\mathbf{z}$, $l \in \Omega^*$ and $0 \le t \le 1$.

M5. A jointly measurable mapping $U: X \times \Omega \to H$, $\pi_\cdot$-nonanticipative in its second variable, and a pair of $\mathcal{F}_\cdot$-adapted random variables $\mathbf{u} \in L^2(P;H)$ and $\mathbf{y}$ ($\Omega$-valued) which satisfy the simultaneous equations

(2.5) $$\begin{cases} \mathbf{y} = \mathbf{u} + \mathbf{w}, \\ \mathbf{u} = U(\mathbf{x}, \mathbf{y}), \end{cases} \qquad P\text{-a.s.}$$

Equivalently $\{(\mathbf{u}_x, \mathbf{y}_x), x \in X\}$ is an $\mathcal{F}_\cdot$-adapted $H \times \Omega$-valued random field with $\mathbf{u}_x \in L^2(P;H)$ $\mu_{\mathbf{x}}$-a.s., and which satisfy

(2.6) $$\begin{cases} \mathbf{y}_x(\theta) = \mathbf{u}_x(\theta) + \mathbf{w}(\theta), \\ \mathbf{u}_x(\theta) = U(x, \mathbf{y}_x(\theta)), \end{cases} \qquad \mu_{\mathbf{x}} \otimes P\text{-a.s.}$$

the connection between (2.5) and (2.6) being $\mathbf{u}(\theta) = \mathbf{u}_x(\theta)|_{x=\mathbf{x}}$.

We now present for later use two facts related to the objects introduced above.

LEMMA 2.1. *For any $h, k \in H$, the function $m(t) := (h, \pi_t k)_H$ is continuous and has bounded variation on $[0,1]$.*

PROOF. The continuity of $m$ follows from that of $t \to \pi_t$. In addition,
$$m(t) = (\pi_t h, \pi_t k)_H = \tfrac{1}{4}(\|\pi_t(h+k)\|^2_H - \|\pi_t(h-k)\|^2_H)$$

so that $m$ has bounded variation, being the difference of two increasing functions. □

LEMMA 2.2. *The random variables of the form (2.3) [thus those of the form (2.4) as well] generate the same $\sigma$-algebra as the one generated by the family of all $(\pi_\cdot, \mathcal{G}_\cdot)$-adapted random variables. This $\sigma$-algebra will be denoted $\mathcal{A}_{\pi_\cdot, \mathcal{G}_\cdot}$.*

PROOF. By density arguments it suffices to check that $\mathbf{0}$ is the only $(\pi_\cdot, \mathcal{G}_\cdot)$-adapted element $\mathbf{u}$ in $L^2(P; H)$ orthogonal to all the random variables of form (2.3). Indeed, for any $s \leq t$ in $[0, 1]$ and $h \in H$ and $a \in L^2(\Theta, \mathcal{F}_s, P)$

$$0 = E(a(\pi_t - \pi_s)h, \mathbf{u})_H = Ea((\pi_t h, u)_H - (\pi_s h, u)_H).$$

This means that $(\pi_t h, u)_H$ is a (continuous) martingale, which in addition has zero bounded variation a.s., by Lemma 2.1. Since it is 0 a.s. for $t = 0$, the same is true for $t = 1$, and since $h \in H$ is arbitrary it follows that $\mathbf{u} = 0$. □

We shall be concerned with the causal and noncausal least mean square estimators

$$(2.7) \qquad \widehat{\mathbf{h}}^{\mathbf{y}} = E(\mathbf{h}|\mathcal{A}_{\pi_\cdot, \mathcal{F}^{\mathbf{y}}_\cdot}) \quad \text{and} \quad \widetilde{\mathbf{h}}^{\mathbf{y}} = E(\mathbf{h}|\mathcal{F}^{\mathbf{y}}_1)$$

of an $H$-valued random variable $\mathbf{h} \in L^2(P; H)$, typically $\mathbf{h} = \mathbf{u}$ or $\mathbf{h} = \mathbf{x}$ (the notation $\mathcal{A}_{\pi_\cdot, \mathcal{F}^{\mathbf{y}}_\cdot}$ was introduced in Lemma 2.2). A central theme of this paper is the relation between their respective associated mean square errors $E|\mathbf{h} - \widehat{\mathbf{h}}^{\mathbf{y}}|_H^2$ and $E|\mathbf{h} - \widetilde{\mathbf{h}}^{\mathbf{y}}|_H^2$ with the mutual information between $\mathbf{x}$ and $\mathbf{y}$, now to be defined.

*Mutual information.* The following definition applies for two general random variables $\mathbf{x}$ and $\mathbf{y}$ defined on a common probability space, the latter taking values in a Polish space so that $\mathbf{y}$'s regular conditional probability measure $\mu_{\mathbf{y}|\mathbf{x}}$ conditioned on $\mathbf{x}$ is well defined. In our case, where $\mathbf{x}$ is given in M2 and $\mathbf{y}$ by the equations (2.5), the key observation is that $\mu_{\mathbf{y}|\mathbf{x}}$ can be expressed in terms of the image measures $\mu_{\mathbf{y}_x}$ of the elements $\mathbf{y}_x, x \in X$, introduced in (2.6):

$$(2.8) \qquad \mu_{\mathbf{y}|\mathbf{x}} = \mu_{\mathbf{y}_x}|_{x=\mathbf{x}} \qquad P\text{-a.s.}$$

DEFINITION 2.3. The mutual information between $\mathbf{x}$ and $\mathbf{y}$ is defined to be

$$(2.9) \qquad I(\mathbf{x}; \mathbf{y}) = \begin{cases} E\left(\log \dfrac{d\mu_{\mathbf{y}|\mathbf{x}}}{d\mu_{\mathbf{y}}}(\mathbf{y})\right), & \text{if } \mu_{\mathbf{y}|\mathbf{x}} \ll \mu_{\mathbf{y}}, P\text{-a.s.} \\ \infty, & \text{otherwise.} \end{cases}$$

Despite (2.9)'s apparent asymmetry, it turns out that $I(\mathbf{x}; \mathbf{y}) = I(\mathbf{y}; \mathbf{x})$. In fact, the identities $\frac{f(y|x)}{f(y)} = \frac{f(x|y)}{f(x)} = \frac{f(x,y)}{f(x)f(y)}$ generalize easily beyond finite dimensions: the following fact is well known and its proof is straightforward.

LEMMA 2.4.

$$\mu_{\mathbf{y}|\mathbf{x}} \ll \mu_{\mathbf{y}}, \quad \mu_{\mathbf{x}}\text{- }a.s. \iff \mu_{\mathbf{x}|\mathbf{y}} \ll \mu_{\mathbf{x}}, \quad \mu_{\mathbf{y}}\text{-}a.s.$$
$$\iff \mu_{\mathbf{x},\mathbf{y}} \ll \mu_{\mathbf{x}} \otimes \mu_{\mathbf{y}}$$

and when one and thus all of these hold, $\frac{d\mu_{\mathbf{y}|\mathbf{x}}}{d\mu_{\mathbf{y}}}(\mathbf{y}) = \frac{d\mu_{\mathbf{x}|\mathbf{y}}}{d\mu_{\mathbf{x}}}(\mathbf{x}) = \frac{d\mu_{\mathbf{x},\mathbf{y}}}{d\mu_{\mathbf{x}} \otimes \mu_{\mathbf{y}}}(\mathbf{x},\mathbf{y})$ $P$-a.s.

Whenever valid (i.e., as long as one does not get $\infty - \infty$) it will be convenient to write

$$(2.10) \qquad I(\mathbf{x};\mathbf{y}) = E \log \frac{d\mu_{\mathbf{y}|\mathbf{x}}}{d\mu_{\mathbf{w}}}(\mathbf{y}) - E \log \frac{d\mu_{\mathbf{y}}}{d\mu_{\mathbf{w}}}(\mathbf{y}),$$

since both terms in the difference can be derived from a generalized Girsanov theorem.

**3. The connection between estimation errors and mutual information.** The main result of this section is the following theorem. It implies in particular that the causal least mean square error does not depend on the resolution of identity which dictates the time structure.

THEOREM 3.1. *Within the setup* M1–M5, *and recalling the notation (2.7)*,

$$(3.1) \qquad I(\mathbf{x},\mathbf{y}) = \tfrac{1}{2} E|\mathbf{u} - \widehat{\mathbf{u}}^{\mathbf{y}}|_H^2$$

*and in the particular case* $\mathbf{y} = \sqrt{\gamma}\mathbf{x} + \mathbf{w}$ *of (2.5)*,

$$(3.2) \qquad I(\mathbf{x},\mathbf{y}) = \frac{\gamma}{2} E|\mathbf{x} - \widehat{\mathbf{x}}^{\mathbf{y}}|_H^2.$$

In the classical case $\Omega = C_0[0,T]$, (3.2) goes back to [1] and the more general case (3.1) in which feedback is allowed was obtained in [5]. The new contribution here is the full extension of (3.1) to the abstract setup. The heart of its proof consists in deriving, in the next subsection, expressions for the Radon–Nikodym derivatives appearing in (2.10) from an abstract version of Girsanov's formula. The theorem's proof will be finalized in Section 3.2.

In this context it is worth stating a recently obtained (for linear observations) connection between the *noncausal* error and mutual information.

THEOREM 3.2 ([3, 13]).  *In the particular case* $\mathbf{y} = \sqrt{\gamma}\mathbf{x} + \mathbf{w}$ *of* (2.5)

$$\frac{dI(\mathbf{x}, \mathbf{y})}{d\gamma} = \frac{1}{2}E|\mathbf{x} - \widetilde{\mathbf{x}}^{\mathbf{y}}|_H^2. \tag{3.3}$$

Theorems 3.1 and 3.2 together yield the following interesting connection between the causal and noncausal errors (cf. [3] as well).

COROLLARY 3.3. *For* $\mathbf{y} = \sqrt{\gamma}\mathbf{x} + \mathbf{w}$ *denote* $\widehat{\varepsilon}^2(\gamma) = E|\mathbf{x} - \widehat{\mathbf{x}}^{\mathbf{y}}|_H^2$ *and* $\widetilde{\varepsilon}^2(\gamma) = E|\mathbf{x} - \widetilde{\mathbf{x}}^{\mathbf{y}}|_H^2$. *Then*

$$\widetilde{\varepsilon}^2(\gamma) = \frac{d(\gamma\widehat{\varepsilon}^2(\gamma))}{d\gamma}, \qquad \text{that is,}$$

$$\widehat{\varepsilon}^2(\gamma) = \frac{\gamma_0 \widehat{\varepsilon}^2(\gamma_0)}{\gamma} + \frac{1}{\gamma}\int_{\gamma_0}^{\gamma} \widetilde{\varepsilon}^2(\beta)\, d\beta \qquad \forall \gamma_0. \tag{3.4}$$

*In addition,* $\widehat{\varepsilon}^2(\frac{1}{\eta})$ *is a concave function of* $\eta$. (We thank an anonymous referee who pointed out an error in an earlier version of this statement, and in its proof.)

PROOF. The identity (3.4) follows directly from (3.2) and (3.3). As for the concavity, denote $h(\eta) = \widehat{\varepsilon}^2(\frac{1}{\eta}) = 2\eta I(\mathbf{x}, \mathbf{y})$. Then

$$h'(\eta) = 2I(\mathbf{x}, \mathbf{y}) + 2\eta\left(-\frac{1}{\eta^2}\right)\frac{dI(\mathbf{x}, \mathbf{y})}{d\gamma} = 2I(\mathbf{x}, \mathbf{y}) - \frac{\widetilde{\varepsilon}^2(1/\eta)}{\eta} \quad \text{and}$$

$$h''(\eta) = \left(-\frac{1}{\eta^2}\right)\frac{dI(\mathbf{x}, \mathbf{y})}{d\gamma} - \frac{1}{\eta}\frac{d}{d\eta}\left(\widetilde{\varepsilon}^2\left(\frac{1}{\eta}\right)\right) + \frac{\widetilde{\varepsilon}^2(1/\eta)}{\eta^2} = -\frac{1}{\eta}\frac{d}{d\eta}\left(\widetilde{\varepsilon}^2\left(\frac{1}{\eta}\right)\right) \leq 0$$

since $\widetilde{\varepsilon}^2(\gamma)$ is clearly a nonincreasing function of $\gamma$.  □

REMARK 3.4.  Viewing $\widehat{\varepsilon}^2$ as a function of $\frac{1}{\eta}$ is equivalent to considering the equally natural model $\mathbf{y} = \mathbf{x} + \sqrt{\eta}\mathbf{w}$ instead of $\mathbf{y} = \sqrt{\gamma}\mathbf{x} + \mathbf{w}$.

3.1. *Girsanov theorem and Radon–Nikodym derivatives on* $\Omega$. Throughout this subsection, $\{\mathcal{G}_t, 0 \leq t \leq 1\}$ will be a generic subfiltration of $\{\mathcal{F}_t, 0 \leq t \leq 1\}$ typically $\mathcal{F}_\cdot^{\mathbf{w}}$ or $\mathcal{F}_\cdot^{\mathbf{y}}$ as defined in (2.2). First, recall the standard Girsanov theorem, in which $\Omega$ is the classical Wiener space $C_0([0, 1])$.

PROPOSITION 3.5.  *Let* $\{b_t, 0 \leq t \leq 1\}$ *be a standard* $\mathcal{G}$*.-Brownian motion,* $\{a_t, 0 \leq t \leq 1\}$ *an* $\mathcal{G}$*.-adapted stochastic process with* $\dot{a}_\cdot \in L^2(0, 1)$ *a.s. and* $y_t = a_t + b_t, 0 \leq t \leq 1$. *Denote*

$$\Lambda_a = \exp\left(-\int_0^1 \dot{a}_t\, db_t - \tfrac{1}{2}\int_0^1 \dot{a}_t^2\, dt\right). \tag{3.5}$$

If $E\Lambda_a = 1$ then $\{y_t, 0 \leq t \leq 1\}$ is a standard $\mathcal{G}.$-Brownian motion. Equivalently (the Jacobi change of variable formula) $E\Lambda_a F(y.) = EF(b.)$, for any $F \in C_b(C_0[0,1])$.

In the context of an abstract Wiener space Itô's integral is defined along the same lines as in the classical case. We now proceed to summarize its construction and refer the reader to [10], Section 2.6, for a more detailed account.

DEFINITION 3.6. An $\Omega$-valued random variable $\mathbf{v}$ is said to be a $\mathcal{G}.$-abstract Wiener process if, for all $l \in \Omega^*$, $M_t^{\mathbf{v},l} := {}_\Omega\langle \mathbf{v}, \pi_t l\rangle_{\Omega^*}$ is a zero mean (and necessarily continuous, Gaussian) $\mathcal{G}.$-martingale with quadratic variation $\langle M^{\mathbf{v},l}\rangle_t = |\pi_t l|_H^2, 0 \leq t \leq 1$.

Note that $\mathbf{w}$ itself is an $\mathcal{F}_\cdot^{\mathbf{w}}$-abstract Wiener process.

Any $\mathcal{G}.$-abstract Wiener process $\mathbf{v}$ generates its associated zero mean Gaussian random field $\{\delta_{\mathbf{v}} l := {}_\Omega\langle \mathbf{v}, l\rangle_{\Omega^*}, l \in \Omega^*\}$ with covariance structure $E\delta_{\mathbf{v}} l_1 \delta_{\mathbf{v}} l_2 = \frac{1}{4}(E_\Omega\langle\mathbf{v}, l_1+l_2\rangle_{\Omega^*}^2 - E_\Omega\langle\mathbf{v}, l_1-l_2\rangle_{\Omega^*}^2) = \frac{1}{4}(|l_1+l_2|_H^2 - |l_1-l_2|_H^2) = (l_1, l_2)_H$ which can thus be extended by density to an $H$-indexed zero mean isonormal Gaussian field $\{\delta_{\mathbf{v}} h, h \in H\}$.

As defined in (i), (ii) and (iii) below, the integrator of Itô's integral will be a $\mathcal{G}.$-abstract Wiener processes $\mathbf{v}$, and in (iv) its "semimartingale" extension. The integrands, now to be defined, will be $(\mathcal{G}., \pi.)$-adapted $H$-valued random variables (eventually all of them).

(i) For $\mathbf{h}$ *simple* as in (2.3), that is, $\mathbf{h} = \sum_{k=0}^{n-1} a_k h_k$ with $0 = t_0 < \cdots < t_n = 1$,
$a_k \in L^2(\Theta, \mathcal{G}_{t_k}, P)$ and $h_k \in (\pi_{t_{k+1}} - \pi_{t_k})(H)$, define $\delta_{\mathbf{v}} \mathbf{h} = \sum_{k=1}^{n} a_k \delta_{\mathbf{v}} h_k$.
For any such $\mathbf{h}$

(3.6) $$E\delta_{\mathbf{v}}\mathbf{h} = 0, \qquad E(\delta_{\mathbf{v}}\mathbf{h})^2 = E|\mathbf{h}|_H^2.$$

(ii) By (3.6) $\delta_{\mathbf{v}}$ can be isometrically extended to the closure in $L^2(P; H)$ of the simple random variables, which turns out to be the set of $(\mathcal{G}., \pi.)$-adapted elements of $L^2(P; H)$. This extension satisfies (3.6) as well.
(iii) For any $(\mathcal{G}., \pi.)$-adapted $H$-valued random variable $\mathbf{h}$, the sequence of $\mathcal{G}.$-stopping times $\tau_n = \inf\{t \in [0,1] \text{ s.t } |\pi_t \mathbf{h}|_H^2 \geq n\}$ ($\inf \varnothing = 1$) increases to 1 as $n \to \infty$, and $\delta_{\mathbf{v}}\mathbf{h} := \lim_{n \to \infty} \delta_{\mathbf{v}} \pi_{\tau_n} \mathbf{h}$ exists almost surely.
(iv) Whenever $\mathbf{z} = \mathbf{v} + \mathbf{u}$, where $\mathbf{u}$ is an $H$-valued random variable, define $\delta_{\mathbf{z}}\mathbf{h} := \delta_{\mathbf{v}}\mathbf{h} + (\mathbf{u}, \mathbf{h})_H$ for any $H$-valued random variable $\mathbf{h}$ which is $(\mathcal{G}., \pi.)$-adapted. This definition is independent of $\mathbf{z}$'s representation.

In abstract Wiener space, (3.5) becomes, for any $(\mathcal{G}., \pi.)$-adapted $\mathbf{h}$,

(3.7) $$\Lambda_{\mathbf{h}} := \exp(-\delta_{\mathbf{v}}\mathbf{h} - \tfrac{1}{2}|\mathbf{h}|_H^2) = \exp(-\delta_{\mathbf{y}}\mathbf{h} + \tfrac{1}{2}|\mathbf{h}|_H^2).$$

PROPOSITION 3.7 ([10]). *Let $\mathbf{v}$ be a $\mathcal{G}.$-abstract Wiener process, $\mathbf{h}$ a $(\mathcal{G}.,\pi.)$-adapted $H$-valued random variable, and $\mathbf{y} = \mathbf{h} + \mathbf{v}$. If $E\Lambda_{\mathbf{h}} = 1$ then $\mathbf{y}$ is a $\mathcal{G}.$ - abstract Wiener process on $(\Theta, \mathcal{F}, \Lambda_{\mathbf{h}} P)$. In particular*

$$(3.8) \qquad E\Lambda_{\mathbf{h}}\varphi(\mathbf{y}) = E\varphi(\mathbf{v}) \qquad \forall \varphi \in C_b(\Omega).$$

*Moreover, $\mathbf{y}$'s and $\mathbf{v}$'s image measures $\mu_{\mathbf{y}}$ and $\mu_{\mathbf{v}}$ are mutually absolutely continuous, and*

$$(3.9) \qquad \frac{d\mu_{\mathbf{y}}}{d\mu_{\mathbf{v}}}(\mathbf{y}) = \frac{1}{E(\Lambda_{\mathbf{h}}|\mathcal{F}_1^{\mathbf{y}})}, \qquad P\text{-a.s.}$$

*[i.e., $\frac{d\mu_{\mathbf{y}}}{d\mu_{\mathbf{v}}}(\omega) = \frac{1}{\lambda_{\mathbf{h}}(\omega)}, \mu_{\mathbf{y}}$-a.s., where $E(\Lambda_{\mathbf{h}}|\mathcal{F}_1^{\mathbf{y}}) = \lambda_{\mathbf{h}}(\mathbf{y})$.]*

PROOF. The Girsanov statement (3.8) is a straightforward generalization of the classical Girsanov theorem (Proposition 3.5), a proof of which can be found in [10], Theorem 2.6.3. From (3.8) it follows for all $\varphi \in C_b(\Omega)$ that $\int_\Omega \varphi(\omega)\mu_{\mathbf{v}}(d\omega) = E\Lambda_{\mathbf{h}}\varphi(\mathbf{y}) = E\lambda_{\mathbf{h}}(\mathbf{y})\varphi(\mathbf{y}) = \int_\Omega \lambda_{\mathbf{h}}(\omega)\varphi(\omega)\mu_{\mathbf{y}}(d\omega)$, and thus $\mu_{\mathbf{v}} \ll \mu_{\mathbf{y}}$ with $\frac{d\mu_{\mathbf{v}}}{d\mu_{\mathbf{y}}} = \lambda_{\mathbf{h}}$, $\mu_{\mathbf{y}}$-a.s.

Moreover, since $\Lambda$ is strictly positive $P$-a.s., so is $\lambda_{\mathbf{h}}, \mu_{\mathbf{y}}$-a.s., and thus $\mu_{\mathbf{v}}$-a.s. as well. This means that $\mu_{\mathbf{y}} \sim \mu_{\mathbf{v}}$ and $\frac{d\mu_{\mathbf{y}}}{d\mu_{\mathbf{v}}} = \frac{1}{\lambda_{\mathbf{h}}}$ $\mu_{\mathbf{y}}$-a.s. □

Although the assumption $E\Lambda_{\mathbf{h}} = 1$ in Proposition 3.7 holds under weaker Novikov-type requirements, the following stronger sufficient condition will suit our needs.

LEMMA 3.8. *Given a $\mathcal{G}.$-abstract Wiener process $\mathbf{v}$, if $\mathbf{h} \in L^\infty(P; H)$ is a $(\mathcal{G}.,\pi.)$-adapted, then $\{\Lambda_{\pi_t\mathbf{h}}, 0 \leq t \leq 1\}$ is a $\mathcal{G}.$-martingale. In particular $E\Lambda_{\mathbf{h}} = 1$.*

PROOF. Assume first that $\mathbf{h} = \sum_{k=1}^n a_k h_k$ is simple, and note that $|a_k| \leq M$ a.s. for some $M < \infty$ and $k = 1, \ldots, n$, and that $E\Lambda_{\pi_t\mathbf{h}} \leq Ee^{M \sum_{k=1}^n |\delta_{\mathbf{v}} h_k|} < \infty$, for any $0 \leq t \leq 1$ since $\delta_{\mathbf{v}} h_1, \ldots, \delta_{\mathbf{v}} h_n$ are Gaussian and independent. To show that $E(\Lambda_{\pi_t\mathbf{h}}|\mathcal{G}_s) = \Lambda_{\pi_s\mathbf{h}}$ for all $s < t$ we may assume without loss of generality that $s = t_{m-1}$ and $t = t_m$ for some $m$. In this case $\Lambda_{\pi_{t_m}\mathbf{h}} = e^{-\sum_{k=1}^m (a_k \delta_{\mathbf{v}} h_k - (1/2) a_k^2 |h_k|_H^2)}$ and

$$\begin{aligned} &E(\Lambda_{\pi_{t_m}\mathbf{h}}|\mathcal{G}_{t_{m-1}}) \\ &= e^{-\sum_{k=1}^{m-1}(a_k \delta_{\mathbf{v}} h_k - (1/2)a_k^2|h_k|_H^2)} E(e^{-a_m \delta_{\mathbf{v}} h_m}|\mathcal{G}_{t_{m-1}}) e^{-(1/2)a_m^2|h_m|_H^2} \\ &= \Lambda_{\pi_{t_{m-1}}\mathbf{h}} \end{aligned}$$

since $a_m$ is $\mathcal{G}_{t_{m-1}}$-measurable and $\delta_{\mathbf{v}} h_m \sim N(0, |\mathbf{h}_m|_H^2)$ is independent of $\mathcal{G}_{t_{m-1}}$.

If $\mathbf{h}$ is $(\mathcal{G}_{\cdot},\pi_{\cdot})$-adapted and $|\mathbf{h}|_H \leq M < \infty$ a.s., let $\mathbf{h}_n$ be a sequence of simple adapted $H$-valued random variables such that $\mathbf{h}_n \to \mathbf{h}$ in $L^2(\theta, \mathcal{F}, P; H)$ as $n \to \infty$. Then $E(\Lambda_{\pi_t \mathbf{h}_n} | \mathcal{G}_s) = \Lambda_{\pi_s \mathbf{h}_n}$ for any $n \in \mathbb{N}$ and $s < t$. Clearly $\Lambda_{\pi_r \mathbf{h}_n} \to \Lambda_{\pi_r \mathbf{h}}$ in probability as $n \to \infty$, for $r = s$ and $r = t$. Since

$$E\Lambda_{\pi_t \mathbf{h}_n}^2 = E\Lambda_{\pi_t 2\mathbf{h}_n} e^{|\pi_t \mathbf{h}_n|_H^2} \leq e^{M^2} E\Lambda_{\pi_t 2\mathbf{h}_n} = e^{M^2}$$

the conditional expectation converges as well and thus $E(\Lambda_{\pi_t \mathbf{h}} | \mathcal{G}_s) = \Lambda_{\pi_s \mathbf{h}}$.
□

PROPOSITION 3.9. *If, in Proposition 3.7, $\mathcal{G}_{\cdot}$ can be taken to be $\mathcal{F}_{\cdot}^{\mathbf{y}}$ (i.e. $\mathbf{v}$ is an $\mathcal{F}_{\cdot}^{\mathbf{y}}$-abstract Wiener process), then $\mu_{\mathbf{y}}$ and $\mu_{\mathbf{v}}$ are mutually absolutely continuous, with*

$$(3.10) \quad \frac{d\mu_{\mathbf{y}}}{d\mu_{\mathbf{v}}}(\mathbf{y}) = \Lambda_{\mathbf{h}}^{-1} \qquad (= e^{\delta_{\mathbf{v}} \mathbf{h} + 1/2|\mathbf{h}|_H^2} = e^{\delta_{\mathbf{y}} \mathbf{h} - 1/2|\mathbf{h}|_H^2}), \qquad P\text{-a.s.}$$

[In this case $\Lambda_{\mathbf{h}}$ is $\mathcal{F}_1^{\mathbf{y}}$-measurable. The point here is that $\frac{d\mu_{\mathbf{y}}}{d\mu_{\mathbf{v}}}(\mathbf{y}) = \Lambda_{\mathbf{h}}^{-1}$, as in (3.9), without requiring $E\Lambda_{\mathbf{h}} = 1$ a priori.] The following proof is essentially taken from [10] Theorem 2.4.2 (where $\mathbf{y}$ is referred to as an *indirect shift* of $\mathbf{v}$) and adapted here to the abstract Wiener space setup.

PROOF. Define $\tau_n = \inf\{t \in [0,1] \text{ s.t } |\pi_t \mathbf{h}|_H \geq n\}$ ($\inf \varnothing = 1$) and let $\mathbf{y}_n = \mathbf{h}_n + \mathbf{v}$ with $\mathbf{h}_n = \pi_{\tau_n} \mathbf{h}$. Since $|\mathbf{h}_n|_H \leq n$ a.s., Lemma 3.8 guarantees that $E\Lambda_{\mathbf{h}_n} = 1$ so that it follows from Proposition 3.7 that $\mu_{\mathbf{y}_n} \sim \mu_{\mathbf{v}}$ and $\frac{d\mu_{\mathbf{y}_n}}{d\mu_{\mathbf{v}}}(\mathbf{y}_n) = \Lambda_{\mathbf{h}_n}^{-1}$ a.s., since $\Lambda_{\mathbf{h}_n}$ itself is $\mathcal{F}_1^{\mathbf{y}}$-measurable, and also $\frac{d\mu_{\mathbf{v}}}{d\mu_{\mathbf{y}_n}}(\mathbf{y}_n) = \Lambda_{\mathbf{h}_n}$. Thus, for any $\varphi \in C_b(\Omega)$,

$$\int_\Theta \varphi \circ \mathbf{v}\, dP = \int_\Omega \varphi\, d\mu_{\mathbf{v}} = \int_\Omega \varphi \frac{d\mu_{\mathbf{v}}}{d\mu_{\mathbf{y}_n}}\, d\mu_{\mathbf{y}_n} = \int_\Theta \varphi \circ \mathbf{y}_n \Lambda_{\mathbf{h}_n}\, dP$$

$$\xrightarrow[n\to\infty]{} \int_\Theta \varphi \circ \mathbf{y} \Lambda_{\mathbf{h}}\, dP$$

since $\Lambda_{\mathbf{h}_n} \to \Lambda_{\mathbf{h}}$ a.s., and thus by Scheffé's lemma $\Lambda_{\mathbf{h}_n} dP \to \Lambda_{\mathbf{h}} dP$ in total variation. This means that $\mu_{\mathbf{v}} \ll \mu_{\mathbf{y}}$ and $\frac{d\mu_{\mathbf{v}}}{d\mu_{\mathbf{y}}}(\mathbf{y}) = \Lambda_{\mathbf{h}}$, and since $\Lambda_{\mathbf{h}} > 0$ a.s., the reverse is true as well, namely $\mu_{\mathbf{y}} \ll \mu_{\mathbf{v}}$ and $\frac{d\mu_{\mathbf{y}}}{d\mu_{\mathbf{v}}}(\mathbf{y}) = \Lambda_{\mathbf{h}}^{-1}$. □

The assumption in Proposition 3.9 that $\mathbf{v}$ is an $\mathcal{F}_{\cdot}^{\mathbf{y}}$-abstract Wiener process suggests that $\mathbf{y} = \mathbf{h} + \mathbf{v}$ should be interpreted as a nonanticipative feedback model, and can thus be expected to hold in the case (2.6):

PROPOSITION 3.10. *Assume the setup M1–M5 in Section 2. Then for $\mu_{\mathbf{x}}$ almost every $x \in X$, $\mathbf{w}$ in (2.6) is an $\mathcal{F}_{\cdot}^{\mathbf{y}_x}$-abstract Wiener process, $\mu_{\mathbf{y}_x} \ll \mu_{\mathbf{w}}$*

*and*

$$\frac{d\mu_{\mathbf{y}_x}}{d\mu_{\mathbf{w}}}(\mathbf{y}_x) = \exp\left(\delta_{\mathbf{w}}\mathbf{u}_x + \frac{1}{2}|\mathbf{u}_x|_H^2\right), \qquad P\text{-}a.s. \tag{3.11}$$

*For the model (2.5),*

$$\frac{d\mu_{\mathbf{y}|\mathbf{x}}}{d\mu_{\mathbf{w}}}(\mathbf{y}) = \exp\left(\delta_{\mathbf{w}}\mathbf{u} + \frac{1}{2}|\mathbf{u}|_H^2\right), \qquad P\text{-}a.s., \tag{3.12}$$

*and in particular*

$$E\log\frac{d\mu_{\mathbf{y}|\mathbf{x}}}{d\mu_{\mathbf{w}}}(\mathbf{y}) = \frac{1}{2}E|\mathbf{u}|_H^2. \tag{3.13}$$

PROOF. Recall that $\mathbf{w}$ is an $\mathcal{F}_\cdot$-abstract Wiener process (cf. Definition 3.6). On the other hand, from (2.6) and bearing in mind that the mapping $U(x,\cdot)$ is nonanticipative, we conclude that $\mathcal{F}_t^{\mathbf{w}} \subset \mathcal{F}_t^{\mathbf{y}_x} \subset \mathcal{F}_t$ for all $0 \leq t \leq 1$, so that $M_t^{\mathbf{w},l} = {}_\Omega\langle\mathbf{w},\pi_t l\rangle_{\Omega^*}$ is not only an $\mathcal{F}_\cdot$-martingale for each $l \in \Omega^*$ but also an $\mathcal{F}_\cdot^{\mathbf{y}_x}$-martingale, and with the same quadratic variation. In other words, $\mathbf{w}$ is an $\mathcal{F}_\cdot^{\mathbf{y}_x}$-abstract Wiener process.

Thus Proposition 3.9 applies to $\mathbf{y}_x = \mathbf{u}_x + \mathbf{w}$ with $\mathbf{h} = \mathbf{u}_x$ and $\mathbf{v} = \mathbf{w}$ [$\mathbf{u}_x$ is indeed $\mathcal{F}^{\mathbf{y}_x}$-adapted, again by (2.6) and $U$'s nonanticipativity], and (3.11) follows.

As for (3.12) we first claim that $\frac{d\mu_{\mathbf{y}|\mathbf{x}}}{d\mu_{\mathbf{w}}}(\mathbf{y}) = \frac{d\mu_{\mathbf{y}_x}}{d}\mu_{\mathbf{w}}(\mathbf{y}_x)|_{x=\mathbf{x}}$. Indeed, for any $\psi \in C_b(\Omega)$,

$$E(\psi(\mathbf{y})|\mathbf{x}) = E(\psi(\mathbf{u_x} + \mathbf{w})|\mathbf{x}) = E\psi(\mathbf{u}_x + \mathbf{w})|_{x=\mathbf{x}} = E\psi(\mathbf{y}_x)|_{x=\mathbf{x}}$$

$$= E\left(\frac{d\mu_{\mathbf{y}_x}}{d\mu_{\mathbf{w}}}(\mathbf{w})\psi(\mathbf{w})\right)\Big|_{x=\mathbf{x}} = E\left(\frac{d\mu_{\mathbf{y}_x}}{d\mu_{\mathbf{w}}}(\mathbf{w})\Big|_{x=\mathbf{x}}\psi(\mathbf{w})\right)$$

(where the independence of $\mathbf{x}$ and $\mathbf{w}$ was used in the second and last equalities), from which it follows that $\mu_{\mathbf{y}|\mathbf{x}} \ll \mu_{\mathbf{w}}, \mu_{\mathbf{x}}$-a.s., and thus $\mu_{\mathbf{y}} \ll \mu_{\mathbf{w}}$, and moreover $\frac{d\mu_{\mathbf{y}|\mathbf{x}}}{d\mu_{\mathbf{w}}}(\mathbf{w}) = \frac{d\mu_{\mathbf{y}_x}}{d\mu_{\mathbf{w}}}(\mathbf{w})|_{x=\mathbf{x}}$. By virtue of the absolute continuity itself,

$$\frac{d\mu_{\mathbf{y}|\mathbf{x}}}{d\mu_{\mathbf{w}}}(\mathbf{y}) = \frac{d\mu_{\mathbf{y}_x}}{d\mu_{\mathbf{w}}}(\mathbf{y})\Big|_{x=\mathbf{x}} = \frac{d\mu_{\mathbf{y}_x}}{d\mu_{\mathbf{w}}}(\mathbf{y}_x)\Big|_{x=\mathbf{x}}$$

as claimed. Combining this with (3.11), and recalling that $\mathbf{u} = \mathbf{u}_x|_{x=\mathbf{x}}$, we obtain

$$\frac{d\mu_{\mathbf{y}|\mathbf{x}}}{d\mu_{\mathbf{w}}}(\mathbf{y}) = \exp\left((\delta_{\mathbf{w}}\mathbf{u}_x)|_{x=\mathbf{x}} + \frac{1}{2}|\mathbf{u}|_H^2\right), \qquad P\text{-a.s.}$$

Note, from the definition of Itô's integral $\delta_{\mathbf{w}}$ and the independence of $\mathbf{w}$ and $\mathbf{x}$, that $(\delta_{\mathbf{w}}\mathbf{u}_x)|_{x=\mathbf{x}} = \delta_{\mathbf{w}}\mathbf{u}$. Thus $E\log\frac{d\mu_{\mathbf{y}|\mathbf{x}}}{d\mu_{\mathbf{w}}}(\mathbf{y}) = E\delta_{\mathbf{w}}\mathbf{u} + \frac{1}{2}E|\mathbf{u}|_H^2 = \frac{1}{2}E|\mathbf{u}|_H^2$ by (3.6), since $|\mathbf{u}|_H$ was assumed to have finite second moment. □

Having found an expression for (2.10)s first term based on the representation (2.6), the starting point for the second term is necessarily (2.5). However, in order to be able to apply Proposition 3.9 in this case (**w** is no longer an $\mathcal{F}^{\mathbf{y}}$-abstract Wiener process) it is necessary to replace (2.5) by **y**'s equivalent innovation representation.

LEMMA 3.11. $\mathbf{n} := \mathbf{y} - \widehat{\mathbf{u}}^{\mathbf{y}} = (\mathbf{u} - \widehat{\mathbf{u}}^{\mathbf{y}}) + \mathbf{w}$ is an $\mathcal{F}^{\mathbf{y}}_{\cdot}$-abstract Wiener process.

PROOF. Let $l \in \Omega^*$. We need to show that $M_t^{\mathbf{n},l} := {}_\Omega\langle \mathbf{n}, \pi_t l\rangle_{\Omega^*}$ is an $\mathcal{F}^{\mathbf{y}}_{\cdot}$-martingale with quadratic variation $|\pi_t l|_H^2$. Indeed,

$$E(M_t^{\mathbf{n},l} - M_s^{\mathbf{n},l}|\mathcal{F}_s^{\mathbf{y}})$$
$$= E({}_\Omega\langle \mathbf{n}, (\pi_t - \pi_s)l\rangle_{\Omega^*}|\mathcal{F}_s^{\mathbf{y}})$$
$$= E({}_\Omega\langle \mathbf{u} - \widehat{\mathbf{u}}^{\mathbf{y}}, (\pi_t - \pi_s)l\rangle_{\Omega^*}|\mathcal{F}_s^{\mathbf{y}}) + E(E({}_\Omega\langle \mathbf{w}, (\pi_t - \pi_s)l\rangle_{\Omega^*}|\mathcal{F}_s)|\mathcal{F}_s^{\mathbf{y}})$$

We shall show that both terms above equal zero, assuming without loss of generality that $s$ and $t$ are dyadic. The second term is indeed zero since **w** is an $\mathcal{F}_{\cdot}$-abstract Wiener process. For the first term, denote by $\mathcal{A}_{\pi_{\cdot},\mathcal{F}^{\mathbf{y},n}_{\cdot}}$ the $\sigma$-algebra generated by the $H$-valued random variables of the form (2.4) on the partition $\mathcal{P} = \{\frac{k}{2^n}, k=0,\ldots,2^n\}$ of $[0,1]$, and

$$\widehat{\mathbf{u}}^{\mathbf{y},n} = E(\mathbf{u}|\mathcal{A}_{\pi_{\cdot},\mathcal{F}^{\mathbf{y},n}_{\cdot}}) = \sum_{k=0}^{2^n-1} E((\pi_{(k+1)/2^n} - \pi_{k/2^n})\mathbf{u}|\mathcal{F}^{\mathbf{y}}_{k/2^n}).$$

It follows from Lemma 2.2 that $\widehat{\mathbf{u}}^{\mathbf{y},n} \to \widehat{\mathbf{u}}^{\mathbf{y}}$ in $L^2(P;H)$, so that it suffices to show that $E({}_\Omega\langle \mathbf{u} - \widehat{\mathbf{u}}^{\mathbf{y},n}, (\pi_t - \pi_s)l\rangle_{\Omega^*}|\mathcal{F}_s^{\mathbf{y}}) = 0$ for every $n$ large enough. Denote $\mathbf{u}_k = (\pi_{(k+1)/2^n} - \pi_{k/2^n})\mathbf{u}$ and, by the dyadic assumption, $s = \frac{k_0}{2^n}$ and $t = \frac{k_1}{2^n}$. Then

$$E({}_\Omega\langle \mathbf{u} - \widehat{\mathbf{u}}^{\mathbf{y},n}, (\pi_t - \pi_s)l\rangle_{\Omega^*}|\mathcal{F}_s^{\mathbf{y}})$$
$$= \sum_{k=k_0}^{k_1-1} E({}_\Omega\langle \mathbf{u}_k - E(\mathbf{u}_k|\mathcal{F}^{\mathbf{y}}_{k/2^n}), l\rangle_{\Omega^*}|\mathcal{F}^{\mathbf{y}}_{k_0/2^n})$$
$$= \sum_{k=k_0}^{k_1-1} E(E({}_\Omega\langle \mathbf{u}_k - E(\mathbf{u}_k|\mathcal{F}^{\mathbf{y}}_{k/2^n}), l\rangle_{\Omega^*}|\mathcal{F}^{\mathbf{y}}_{k/2^n})|\mathcal{F}^{\mathbf{y}}_{k_0/2^n}) = 0.$$

As for the quadratic variation, note that $M_t^{\mathbf{n},l} = (\mathbf{u} - \widehat{\mathbf{u}}^{\mathbf{y}}, \pi_t l)_H + M_t^{\mathbf{w},l}$. By Lemma 2.1 the first term is almost surely continuous, has bounded variation and thus zero quadratic variation, so that $\langle M^{\mathbf{n},l}\rangle_t = \langle M^{\mathbf{w},l}\rangle_t = |\pi_t l|_H^2$. □

COROLLARY 3.12.

$$E \log \frac{d\mu_\mathbf{y}}{d\mu_\mathbf{w}} = \frac{1}{2} E |\widehat{\mathbf{u}}^\mathbf{y}|_H^2. \tag{3.14}$$

PROOF. We may apply Proposition 3.9 to $\mathbf{y} = \widehat{\mathbf{u}}^\mathbf{y} + \mathbf{n}$ to conclude that

$$\frac{d\mu_\mathbf{y}}{d\mu_\mathbf{w}} = \frac{d\mu_\mathbf{y}}{d\mu_\mathbf{n}} = \exp\left(\delta_\mathbf{n} \widehat{\mathbf{u}}^\mathbf{y} + \frac{1}{2}|\widehat{\mathbf{u}}^\mathbf{y}|_H^2\right). \tag{3.15}$$

(Indeed, $\widehat{\mathbf{u}}^\mathbf{y}$ is clearly $\mathcal{F}_\cdot^\mathbf{y}$-adapted and $\mathbf{n}$ is an $\mathcal{F}_\cdot^\mathbf{y}$-abstract Wiener process by Lemma 3.11.) Since $E(\widehat{\mathbf{u}}^\mathbf{y})^2 \leq E\mathbf{u}^2 < \infty$, and thus $E\delta_\mathbf{n}\widehat{\mathbf{u}}^\mathbf{y} = 0$, (3.15) implies (3.14). □

3.2. *Proof of Theorem* 3.1. All that remains to prove (3.1) [and thus (3.2) as well] is to insert (3.13) and (3.14) in (2.10) and thus obtain $I(\mathbf{x}, \mathbf{y}) = \frac{1}{2}E|\mathbf{u}|_H^2 - \frac{1}{2}E|\widehat{\mathbf{u}}^\mathbf{y}|_H^2 = \frac{1}{2}E|\mathbf{u} - \widehat{\mathbf{u}}^\mathbf{y}|_H^2$.

**4. Gaussian signals.** Consider the particular case of (2.5)

$$\mathbf{y} = \sqrt{\gamma}\mathbf{x} + \mathbf{w}, \tag{4.1}$$

where $\mathbf{x}$ is assumed to be a zero mean Gaussian $H$-valued random variable with correlation bilinear form $r(h,k) = E(\mathbf{x},h)_H(\mathbf{x},k)_H, h, k \in H$, and associated correlation operator $\mathbf{R}$ on $H$ characterized by $(\mathbf{R}h,k)_H = r(h,k)$ for all $h, k \in H$. The positive constant $\gamma$ is commonly called the signal to noise ratio.

It is well known that $\mathbf{R}$ is nonnegative and of trace class. Its spectrum thus consists of a nonincreasing summable sequence $\{\lambda_i\}_{i=1}^\infty$ of nonnegative eigenvalues with an associated family $\{\varphi_i\}_{i=1}^\infty$ of orthonormal eigenvectors and

$$\mathbf{R} = \sum_{i=1}^\infty \lambda_i \varphi_i \otimes \varphi_i, \quad \text{that is, } r(h,k) = \sum_{i=1}^\infty \lambda_i (\varphi_i, h)_H (\varphi_i, k)_H \quad \forall h, k \in H,$$

which leads immediately to the representation

$$\mathbf{x} = \sum_{i=1}^\infty \sqrt{\lambda_i} \xi_i \varphi_i, \quad \text{in } L^2(P; H) \tag{4.2}$$

where $\{\xi_i = (\mathbf{x}, \varphi_i)_H\}_{i=1}^\infty$ is an i.i.d. $N(0,1)$ sequence.

THEOREM 4.1. *The least causal mean square error of* $\mathbf{x} \sim N(0, \mathbf{R})$ *with* $\mathbf{y}$ *as in (4.1) is given by*

$$\widehat{\varepsilon}^2(\gamma) = E|\mathbf{x} - \widehat{\mathbf{x}}|^2 = \gamma^{-1} \sum_{i=1}^\infty \log(1 + \lambda_i \gamma) = \gamma^{-1} \log \det(I + \gamma \mathbf{R}). \tag{4.3}$$

If **x** *is only assumed to possess a covariance* **R** *(but not necessarily to be Gaussian), then the right-hand side of (4.3) yields the least linear causal mean square error.*

PROOF. Expanding **y** and **w** in the vectors $\{\varphi_i\}$, one obtains $\eta_i = \sqrt{\gamma}\xi_i + \omega_i$ where $\omega_i = (\mathbf{w}, \varphi_i)_H$ and $\eta_i = (\mathbf{y}, \varphi_i)_H$ are independent for all $i$. From the orthogonality one concludes that

$$(4.4) \quad \tilde{\varepsilon}^2(\gamma) = \sum_{i=1}^{\infty} E(\xi_i - E(\xi_i|\eta_i))^2 = \sum_{i=1}^{\infty} \frac{\lambda_i}{1 + \lambda_i \gamma}$$

(where the last equality is a standard one-dimensional calculation). Applying (3.4) with $\gamma_0 = 0$ we obtain (4.3) as claimed:

$$\hat{\varepsilon}^2(\gamma) = \frac{1}{\gamma} \sum_{i=1}^{\infty} \int_0^{\gamma} \frac{\lambda_i}{1 + \lambda_i \gamma'} d\gamma' = \frac{1}{\gamma} \sum_{i=1}^{\infty} \log(1 + \lambda_i \gamma). \quad \square$$

Note that the formulae (4.3) and (4.4) yield the asymptotic expansions in powers of $\gamma$ in terms of the coefficients $s_k = \sum_i \lambda_i^k$ ($s_k^{1/k}$ are known as **R**'s Schatten norms):

$$\hat{\varepsilon}^2(\gamma) \sim \sum_{k=0}^{\infty} (-1)^k \frac{s_k}{k+1} \gamma^k \quad \text{and} \quad \tilde{\varepsilon}^2(\gamma) \sim \sum_{k=0}^{\infty} (-1)^k s_k \gamma^k.$$

It is of course not surprising that $\hat{\varepsilon}^2(\gamma) \sim \tilde{\varepsilon}^2(\gamma) \sim s_1 = E|\mathbf{x}|_H^2$ as $\gamma \to 0$. A more interesting consequence of these expansions is

COROLLARY 4.2.

$$\lim_{\gamma \to 0} \frac{E|\mathbf{x}|_H^2 - \hat{\varepsilon}^2(\gamma)}{E|\mathbf{x}|_H^2 - \tilde{\varepsilon}^2(\gamma)} = 2.$$

In other words, the noncausal error increases to its limit in small signal to noise ratio twice as fast as the causal error, regardless of the correlation operator. This is not necessarily true if **x** is not assumed to be Gaussian.

The last application of Theorem 4.1 concerns the mean square causal estimation error of a stationary multivariate Gaussian process $\{x_t, t \in \mathbb{R}\}$ in additive white noise. The so called Yovits–Jackson formula for this quantity has been obtained in the scalar case under various assumptions and by different analytic methods, as explained in the Introduction. Here it follows in full generality as a straightforward consequence of Theorem 4.1.

PROPOSITION 4.3. *Let $\{x_t, t \in \mathbb{R}\}$ be a stationary zero mean $n$-dimensional Gaussian process with continuous correlation function $R(\tau) := Ex_0 x_\tau^T \in L^1(\mathbb{R};$*

$\mathbb{R}^{n \times n}$) *and spectral density* $S(\xi)$, *and let* $y_t = \sqrt{\gamma} \int_0^t x_s\, ds + w_t$, $t \in \mathbb{R}$, *where* $\{w_t, t \in \mathbb{R}\}$ *is a two-sided standard n-dimensional Brownian motion and* $\gamma > 0$. *Furthermore, denote by* $y_a^b$ *the sigma algebra generated by* $\{y_t - y_s, a \leq s < t \leq b\}$, *for any* $-\infty \leq a < b \leq \infty$. *Then, for any fixed time* $\theta$,

$$(4.5) \qquad E|x_\theta - E(x_\theta | y_{-\infty}^\theta)|^2 = (2\pi\gamma)^{-1} \int_{-\infty}^{\infty} \log \det(I + \gamma S(\xi))\, d\xi.$$

PROOF. On each finite time interval $[0, T]$, this case can be modeled by the classical Wiener space $\Omega = C_0([0, T]; \mathbb{R}^n)$ with $\mathbf{w} = w.$, $\mathbf{x} = \int_0^{\cdot} x_t\, dt$ and $\mathbf{y} = y. = \sqrt{\gamma} \mathbf{x} + \mathbf{w}$. Let $\mathbf{R}_T$ be the Toeplitz integral operator with kernel $R(t-s)$ and spectrum $\{\lambda_i^{(T)}\}_{i=1}^{\infty}$, and $\mathbf{I}_T$ the identity operator, on $L^2([0, T]; \mathbb{R}^n)$. By Theorem 4.1, and in view of the stationarity,

$$(4.6) \quad \frac{1}{T} \int_0^T E|x_\theta - E(x_\theta | y_{\theta-t}^\theta)|^2\, dt = \frac{1}{T} \int_0^T E|x_t - E(x_t | y_0^t)|^2\, dt$$
$$= \frac{1}{\gamma} \left( \frac{1}{T} \sum_{i=1}^{\infty} \log(1 + \gamma \lambda_i^{(T)}) \right)$$

The integrand in the left-hand side converges, as $t \to \infty$, to the left-hand side of (4.5) by standard martingale theory, and thus so does the integral average itself. The convergence of right-hand side is a consequence of a matrix-valued version of the Kac–Murdock–Szegö theorem on $\mathbf{R}_T$'s asymptotic eigenvalue distribution (see [4], Section 4.4 or [2], page 139). Specifically, [6] Theorem 3.2, states (formula (3.2) in [6] contains a typographical error; the integrand there should be $\operatorname{tr} \Phi(K(t))$, as is evident throughout the subsequent proof) that as $T \to \infty$ the term in parenthesis converges to $\frac{1}{2\pi} \int_{-\infty}^{\infty} \log \det(I + \gamma S(\xi))\, d\xi$ which concludes the proof. (The cited theorem was applied to the function $\Phi(z) = \log(1 + \gamma z)$, $z \in [0, E|x_0|^2]$, which is allowed in view of Remark 3.2 in [6].) □

FACULTY OF MATHEMATICS  
TECHNION, ISRAEL INSTITUTE OF TECHNOLOGY  
HAIFA 32000  
ISRAEL  
E-MAIL: emw@tx.technion.ac.il

FACULTY OF ELECTRICAL ENGINEERING  
TECHNION, ISRAEL INSTITUTE OF TECHNOLOGY  
HAIFA 32000  
ISRAEL  
E-MAIL: zakai@ee.technion.ac.il